# Solution of System of Linear Fractional Differential Equations with Modified derivative of Jumarie Type

**Uttam Ghosh[1] , Susmita Sarkar[2] and Shantanu Das [3]**

[1] Department of Mathematics, Nabadwip Vidyasagar College, Nabadwip, Nadia, West Bengal, India;
Email : uttam_math@yahoo.co.in

[2]Department of Applied Mathematics, University of Calcutta, Kolkata, India
Email : susmita62@yahoo.co.in

[3]Scientist H+, Reactor Control Systems Design Section E & I Group B.A.R.C Mumbai India
Senior Research Professor, Dept. of Physics, Jadavpur University Kolkata
Adjunct Professor. DIAT-Pune
UGC Visiting Fellow. Dept of Appl. Mathematics; Univ. of Calcutta
Email : shantanu@barc.gov.in

## Abstract

Solution of fractional differential equations is an emerging area of present day research because such equations arise in various applied fields. In this paper we have developed analytical method to solve the system of fractional differential equations in-terms of Mittag-Leffler function and generalized Sine and Cosine functions, where the fractional derivative operator is of Jumarie type. The use of Jumarie type fractional derivative, which is modified Rieman-Liouvellie fractional derivative, eases the solution to such fractional order systems. The use of this type of Jumarie fractional derivative gives a conjugation with classical methods of solution of system of linear integer order differential equations, by usage of Mittag-Leffler and generalized trigonometric functions. The ease of this method and its conjugation to classical method to solve system of linear fractional differential equation is appealing to researchers in fractional dynamic systems. Here after developing the method, the algorithm is applied in physical system of fractional differential equation. The analytical results obtained are then graphically plotted for several examples for system of linear fractional differential equation.

## 1.0 Introduction

The fractional calculus is a current research topic in applied sciences such as applied mathematics, physics, mathematical biology and engineering. The rule of fractional derivative is not unique till date. The definition of fractional derivative is given by many authors. The commonly used definition is the Riemann-Liouvellie (R-L) definition [1-5]. Other useful definition includes Caputo definition of fractional derivative (1967) [1-5]. Jumarie's left handed modification of R-L fractional derivative is useful to avoid nonzero fractional derivative of a constant functions [7]. Recently in the paper [8] Ghosh *et al* proposed a theory of characterization of non-differentiable points with Jumarie type fractional derivative with right handed modification of R-L fractional derivative. The differential equations in different form of fractional derivatives give different type of solutions [1-5]. Therefore, there is no standard algorithm to solve fractional differential equations. Thus the solution and its interpretation of the fractional differential equations is a rising field of Applied Mathematics. To solve the linear and non-linear differential equations recently used methods are Predictor-Corrector method [9], Adomain decomposition method [2, 10-11], Homotopy Perturbation Method [12] Variational Iteration Method [13], Differential transform method [14]. Recently in [15] Ghosh *et al* developed analytical method for solution of linear fractional differential equations with Jumarie type derivative [7] in terms of Mittag-Leffler functions and generalized sine and cosine functions. This new finding of [15] has been extended in this paper to get analytical solution of system of linear fractional differential equations. In section 1.0 we have defined some important definitions of fractional derivative that is basic Riemann-Liouvellie (RL) fractional derivative, the Caputo fractional derivative, the Jumarie fractional derivative, the Mittag-Leffler function and generalized Sine and Cosine functions. In section 2.0 solution of system of fractional differential equations has been described and in section 3.0 an application of this method to physical system has been discussed.

**1.1 The basic definitions of fractional derivatives and some higher transcendental functions:**

**a) Basic definitions of fractional derivative:**

**i) Riemann- Liouvellie (R-L) definition**

The R-L definition of the left fractional derivative is,

$$ {}_aD_x^{\alpha} f(x) = \frac{1}{\Gamma(m+1-\alpha)}\left(\frac{d}{dx}\right)^{m+1}\int_a^x (x-\tau)^{m-\alpha} f(\tau)d\tau \quad \text{..............(1.1)}$$

Where: $m \le \alpha < m+1$, $m$: is positive integer.

In particular when $0 \le \alpha < 1$ then

$$ {}_aD_x^{\alpha} f(x) = \frac{1}{\Gamma(1-\alpha)}\frac{d}{dx}\int_a^x (x-\tau)^{-\alpha} f(\tau)d\tau \quad \text{..............(1.2)}$$

The definition (1.1) is known as the left R-L definition of the fractional derivative. The corresponding right R-L definition is

$$ {}_xD_b^{\alpha} f(x) = \frac{1}{\Gamma(m+1-\alpha)}\left(-\frac{d}{dx}\right)^{m+1}\int_x^b (\tau-x)^{m-\alpha} f(\tau)d\tau \quad \text{..............(1.3)}$$

where: $m \le \alpha < m+1$.

The derivative of a constant is obtained as non-zero using the above definitions (1.1)-(1.3) which contradicts the classical derivative of the constant, which is zero. In 1967 Prof. M. Caputo proposed a modification of the R-L definition of fractional derivative which can overcome this shortcoming of the R-L definition.

**ii) Caputo definition**

M. Caputo defines the fractional derivative in the following form [6]

$$ {}_a^CD_x^{\alpha} f(x) = \frac{1}{\Gamma(n-\alpha)}\int_a^x (x-\tau)^{n-\alpha-1} f^{(n)}(\tau)d\tau \quad \text{..............(1.4)}$$

where: $n-1 \le \alpha < n$.

In this definition first differentiate $f(x)$, $n-$times then integrate $n-\alpha$ times. The disadvantage of this method is that $f(x)$, must be differentiable $n$-times then the $\alpha$ -th order derivative will exist, where $n-1 \le \alpha < n$ . If the function is non-differentiable then this definition is not applicable. Two main advantages of this method are (i) fractional derivative of a constant is zero (ii) the fractional differential equation of Caputo type has initial conditions of classical derivative type but the R-L type differential equations has initial conditions fractional type i.e. $\lim_{x\to a} {}_aD_x^{\alpha-1}\left[f(x)\right]=b_1$.

This means that a fractional differential equation composed with RL fractional derivatives require concept of fractional initial states, sometimes they are hard to interpret physically [2].

**iii) Modified definitions of fractional derivative :**

To overcome the fractional derivative of a constant, non-zero, another modification of the definition of left R-L type fractional derivative of the function $f(x)$, in the interval $[a,b]$ was proposed by Jumarie [7] in the form, that is following.

$$f_L^{(\alpha)}(x)={}_a^JD_x^\alpha f(x)=\begin{cases}\dfrac{1}{\Gamma(-\alpha)}\displaystyle\int_a^x(x-\tau)^{-\alpha-1}f(\tau)d\tau, & \alpha<0.\\ \dfrac{1}{\Gamma(1-\alpha)}\dfrac{d}{dx}\displaystyle\int_a^x(x-\tau)^{-\alpha}\left(f(\tau)-f(a)\right)d\tau, & 0<\alpha<1\\ \left(f^{(\alpha-m)}(x)\right)^{(m)}, & m\le\alpha<m+1.\end{cases} \quad .......(1.5)$$

We consider that $f(x)-f(a)=0$ for $x<a$. In (1.5), the first expression is just fractional integration; the second line is RL derivative of order $0<\alpha<1$ of offset function that is $f(x)-f(a)$. For $\alpha>1$, we use the third line; that is first we differentiate the offset function with order $0<(\alpha-m)<1$, by the formula of second line, and then apply whole $m$ order differentiation to it. Here we chose integer $m$, just less than the real number $\alpha$; that is $m\le\alpha<m+1$.

The logic of Jumarie fractional derivative is that, we do RL fractional derivative operation on a new function by forming that new function from a given function by offsetting the value of the function at the start point. Here the differentiability requirement as demanded by Caputo definition is not there. Also the fractional derivative of constant function is zero, which is non-zero by RL fractional derivative definition.

We have recently modified the right R-L definition of fractional derivative of the function $f(x)$, in the interval $[a,b]$ in the following form [8],

$$f_R^{(\alpha)}(x)={}_x^JD_b^\alpha f(x)=\begin{cases}-\dfrac{1}{\Gamma(-\alpha)}\displaystyle\int_x^b(\tau-x)^{-\alpha-1}f(\tau)d\tau, & \alpha<0.\\ -\dfrac{1}{\Gamma(1-\alpha)}\dfrac{d}{dx}\displaystyle\int_x^b(\tau-x)^{-\alpha}\left(f(b)-f(\tau)\right)d\tau, & 0<\alpha<1\\ \left(f^{(\alpha-m)}(x)\right)^{(m)}, & m\le\alpha<m+1.\end{cases} \quad .......(1.6)$$

In the same paper [8], we have shown that both the modifications (1.5) and (1.6) give fractional derivatives of non-differentiable points their values are different, at that point, but we get finite values there, of fractional derivatives. Whereas in classical integer order calculus, where we have different values of right and left derivatives at non differentiable points in approach limit from left side or right side, but infinity (or minus infinity) at that point, where function is non-differentiable. But in case of Jumarie fractional derivative and right modified RL fractional derivative [8], there is no approach limit at the non-differentiable points, but a finite value is obtained at that non differentiable point of the function. The difference is that integer order calculus returns infinity or minus infinity at non-differentiable points, where as the Jumarie fractional derivative returns a finite number indicating the character of otherwise non-differentiable points in a function, in left sense or right sense. This has a significant application in characterizing otherwise non-differentiable but continuous points in the function. However, the finite value of the non differentiable point after fractional differentiation depends on the interval length. In the rest of the paper ${}^{J}D^{\upsilon}$ will represent Jumarie fractional derivative.

b) **Mittag-Leffler function and the generalized Sine and Cosine functions**

The Mittag-Leffler function was introduced by the Swedish mathematician Gösta Mittag-Leffler [17-20] in 1903. It is the direct generalization of exponential functions. The one parameter Mittag-Leffler function is defined (in series form) as:

$$E_\alpha(z) \stackrel{\text{def}}{=} \sum_{k=0}^{\infty} \frac{z^k}{\Gamma(1+\alpha k)}, \qquad z \in \mathbb{C}, \qquad \mathrm{Re}(\alpha) > 0$$

In the solutions of FDE we use this series definition in MATLAB plots. One parameter Mittag-Leffler function in relation to few transcendental functions is as follows

$$\begin{aligned} E_1(z) &\stackrel{\text{def}}{=} 1 + \frac{z}{\Gamma(1+1)} + \frac{z^2}{\Gamma(1+2)} + \frac{z^3}{\Gamma(1+3)} + \ldots\ldots\ldots\ldots\ldots\ldots \\ &= 1 + \frac{z}{1!} + \frac{z^2}{2!} + \frac{z^3}{3!} + \ldots\ldots\ldots\ldots\ldots\ldots \\ &= e^z \end{aligned}$$

$$\begin{aligned} E_2(z) &\stackrel{\text{def}}{=} 1 + \frac{z}{\Gamma(1+2)} + \frac{z^2}{\Gamma(1+4)} + \frac{z^3}{\Gamma(1+6)} + \ldots\ldots\ldots\ldots\ldots\ldots \\ &= 1 + \frac{z}{2!} + \frac{z^2}{4!} + \frac{z^3}{6!} + \ldots\ldots\ldots\ldots\ldots\ldots \\ &= \cosh(\sqrt{z}) \end{aligned}$$

The integral representation of the Mittag-Leffler function [17]-[20] is,

$$E_\alpha(z) = \frac{1}{2\pi} \int_C \frac{t^{\alpha-1} e^t}{t^\alpha - z} dt, \qquad z \in \mathbb{C}, \qquad \mathrm{Re}(\alpha) > 0$$

Here the path of the integral $C$ is a loop which starts and ends at $-\infty$ and encloses the circles of disk $|t|\le|z|^{1/\alpha}$ in positive sense : $|\arg(t)|\le\pi$ on $C$ [17]-[20].

The two parameter Mittag-Leffler function (in series form) and its relation with few transcendental functions are as following

$$E_{\alpha,\beta}(z)\overset{\text{def}}{=}\sum_{k=0}^{\infty}\frac{z^k}{\Gamma(\beta+\alpha k)},\qquad z,\beta\in\mathbb{C},\qquad \operatorname{Re}(\alpha)>0$$

$$E_{\alpha,1}(z)=E_{\alpha}(z)=e^z\qquad z\in\mathbb{C},\qquad \operatorname{Re}(\alpha)>0$$

$$E_{1,2}(z)=\frac{e^z-1}{z}$$

$$E_{2,2}(z)=\frac{\sinh(\sqrt{z})}{\sqrt{z}}.$$

The corresponding integral representation [17]-[20] of the two parameter Mittag-Leffler function is,

$$E_{\alpha,\beta}(z)=\frac{1}{2\pi}\int_C\frac{t^{\alpha-\beta}e^t}{t^\alpha-z}dt,\qquad z\in\mathbb{C},\qquad \operatorname{Re}(\alpha)>0$$

where the contour $C$ is already defined, in the above paragraph.

Using the modified definition, of fractional derivative of Jumarrie type, [7-8] we get

$${}^{J}D^{\alpha}\left[1\right]=0,\qquad 0<\alpha<1.$$

The Jumarie fractional derivative of any constant function is zero, unlike a non-zero value of fractional RL derivative of a constant.

We now find Jumarie fractional derivative of Mittag-Leffler function $E_\alpha(at^\alpha)$

$$E_\alpha(at^\alpha)\overset{\text{def}}{=}1+\frac{at^\alpha}{\Gamma(1+\alpha)}+\frac{a^2t^{2\alpha}}{\Gamma(1+2\alpha)}+\frac{a^3t^{3\alpha}}{\Gamma(1+3\alpha)}+\dots\dots\dots\dots\dots\dots$$

Using Jumarie derivative of order $\alpha$, with $0\le\alpha<1$ with start point as $a=0$ for $f(t)=t^{n\alpha}$, [15] that is

$${}_0^JD_t^\alpha\left[t^{n\alpha}\right]=\frac{\Gamma(n\alpha+1)}{\Gamma\left(\alpha(n-1)+1\right)}(t)^{\alpha(n-1)}$$

for $n=1,2,3,\dots$; and also using Jumarie derivative of constant as zero ${}^JD^\alpha\left[1\right]=0$, we get the following very useful identity.

$$
{}_0^J D_t^\alpha \left[ E_\alpha (at^\alpha) \right] = D^\alpha \left( 1 + \frac{at^\alpha}{\Gamma(1+\alpha)} + \frac{a^2 t^{2\alpha}}{\Gamma(1+2\alpha)} + \frac{a^3 t^{3\alpha}}{\Gamma(1+3\alpha)} + \dots\dots\dots\dots\dots\dots \right)
$$

$$
= 0 + \frac{\Gamma(1+\alpha)a}{\Gamma(1)\Gamma(1+\alpha)} + \frac{\Gamma(1+2\alpha)a^2 t^\alpha}{\Gamma(1+2\alpha)\Gamma(1+\alpha)} + \frac{\Gamma(1+3\alpha)a^3 t^{2\alpha}}{\Gamma(1+3\alpha)\Gamma(1+2\alpha)} + \dots\dots\dots\dots
$$

$$
= a\left( 1 + \frac{at^\alpha}{\Gamma(1+\alpha)} + \frac{a^2 t^{2\alpha}}{\Gamma(1+2\alpha)} + \frac{a^3 t^{3\alpha}}{\Gamma(1+3\alpha)} + \dots\dots\dots\dots\dots \right)
$$

$$
= aE_\alpha (at^\alpha)
$$

Thus

$$
{}_0^J D_t^\alpha \left[ E_\alpha (at^\alpha) \right] = aE_\alpha (at^\alpha)
$$

This shows that $AE_\alpha (at^\alpha)$ is a solution is a solution of the fractional differential equation [15]

$$
{}_0^J D^\alpha y = ay
$$

Where $A$ is arbitrary constant.

Therefore

$$
{}_0^J D^\alpha y = ay
$$

with $y(0) = 1$ has solution

$$
y = E_\alpha (at^\alpha)
$$

The fractional Sine and Cosine functions are expressed as following [16],

$$
E_\alpha (ix^\alpha) \overset{\text{def}}{=} \cos_\alpha (x^\alpha) + i \sin_\alpha (x^\alpha)
$$

$$
\cos_\alpha (x^\alpha) \overset{\text{def}}{=} \sum_{k=0}^{\infty} (-1)^k \frac{x^{2k\alpha}}{\Gamma(1+2\alpha k)}
$$

$$
\sin_\alpha (x^\alpha) \overset{\text{def}}{=} \sum_{k=0}^{\infty} (-1)^k \frac{x^{(2k+1)\alpha}}{\Gamma\left(1+(1+2k)\alpha\right)}
$$

The above series form of fractional Sine and Cosine are used to plots, in solutions. It can be easily shown that [15], [16]

$$
{}^J D^\alpha \left[ \sin_\alpha (x^\alpha) \right] = \cos_\alpha (x^\alpha) \qquad {}^J D^\alpha \left[ \cos_\alpha (x^\alpha) \right] = -\sin_\alpha (x^\alpha)
$$

This has been proved by the following term by term differentiation. The series presentation of $\cos_\alpha (x^\alpha)$ is [15],

$$\cos_\alpha(x^\alpha) = 1 - \frac{x^{2\alpha}}{\Gamma(1+2\alpha)} + \frac{x^{4\alpha}}{\Gamma(1+4\alpha)} - \frac{x^{6\alpha}}{\Gamma(1+6\alpha)} + .....$$

Taking its term by term Jumarie fractional derivative of order $\alpha$ we get,

$$\begin{aligned} {}_0^J D_x^\alpha\left[\cos_\alpha(x^\alpha)\right] &= 0 - \frac{\Gamma(1+2\alpha)x^{2\alpha-\alpha}}{\Gamma(1+2\alpha)\Gamma(1+\alpha)} + \frac{\Gamma(1+4\alpha)x^{4\alpha-\alpha}}{\Gamma(1+4\alpha)\Gamma(1+3\alpha)} - \frac{\Gamma(1+6\alpha)x^{6\alpha-\alpha}}{\Gamma(1+6\alpha)\Gamma(1+5\alpha)} + ..... \\ &= -\left[\frac{x^\alpha}{\Gamma(1+\alpha)} - \frac{x^{3\alpha}}{\Gamma(1+3\alpha)} + .......\right] \\ &= -\sin_\alpha(x^\alpha) \end{aligned}$$

Similarly we can get the expression for fractional derivative of Jumarie type of order $\alpha$ for $\sin_\alpha(x^\alpha)$.

## 2.0 System of linear fractional differential equations

Before considering the system of fractional differential equations we state the results [15] which arises in solving a single linear fractional differential equations composed by Jumarie derivative; we will be using the following theorems.

**Theorem 1:** The fractional differential equation

$$\left({}^J D^\alpha - a\right)\left({}^J D^\alpha - b\right) y(t) = 0$$

has solution of the form

$$y = AE_\alpha(at^\alpha) + BE_\alpha(bt^\alpha)$$

where *A* and *B* are constants.

**Theorem 2:** The fractional differential equation

$${}^J D^{2\alpha}\left[y\right] - 2a\left({}^J D^\alpha\left[y\right]\right) + a^2 y = 0$$

has solution of the form

$$y = (At^\alpha + B)E_\alpha(at^\alpha)$$

where *A* and *B* are constants.

**Theorem 3:** Solution of the fractional differential equation

$$ {}^{J}D^{2\alpha}\left[y\right]-2a\left({}^{J}D^{\alpha}\left[y\right]\right)+(a^{2}+b^{\alpha})y=0 $$

is of the form

$$ y=\left(E_{\alpha}(at^{\alpha})\right)\left(A\cos_{\alpha}(bt^{\alpha})+B\sin_{\alpha}(bt^{\alpha})\right) $$

where $A$ and $B$ are constants.

Consider the system of linear fractional differential equations

$$ \left.\begin{aligned} {}^{J}D^{\alpha}\left[x\right]=ax+by \\ {}^{J}D^{\alpha}\left[y\right]=cx+dy \end{aligned}\right\} \qquad (2.1) $$

Here $a$, $b$, $c$ and $d$ are constants, the operator ${}^{J}D^{\alpha}$ is the Jumarie fractional derivative operator, call it for convenience ${}^{J}D^{\alpha}\equiv\frac{d^{\alpha}}{dt^{\alpha}}$, and $x$ and $y$ are functions of $t$. In matrix form we write the (2.1) in following way

$$ \begin{bmatrix} {}^{J}D^{\alpha}x \\ {}^{J}D^{\alpha}y \end{bmatrix}=\begin{bmatrix} a & b \\ c & d \end{bmatrix}\begin{bmatrix} x \\ y \end{bmatrix} \quad \text{also} \quad \frac{d^{\alpha}}{dt^{\alpha}}\begin{bmatrix} x \\ y \end{bmatrix}=\begin{bmatrix} a & b \\ c & d \end{bmatrix}\begin{bmatrix} x \\ y \end{bmatrix} $$

The above system (2.1) can be written in the following form

$$ \left.\begin{aligned} \left({}^{J}D^{\alpha}-a\right)x-by=0 \\ -cx+\left({}^{J}D^{\alpha}-d\right)y=0 \end{aligned}\right\} \qquad (2.2) $$

Operating $({}^{J}D^{\alpha}-a)$ on both sides of the second equation of (2.2) we get the following steps.

$$ -c({}^{J}D^{\alpha}-a)x+({}^{J}D^{\alpha}-a)({}^{J}D^{\alpha}-d)y=0 $$
$$ -cby+({}^{J}D^{\alpha}-a)({}^{J}D^{\alpha}-d)y=0 $$

$$ {}^{J}D^{2\alpha}\left[y\right]-(a+d)\left({}^{J}D^{\alpha}\left[y\right]\right)+(ad-bc)y=0 \qquad (2.3) $$

Equation (2.3) is a linear fractional order differential equation (with order $2\alpha$), with Jumarie derivative operator.

Let

$$ a+d=\lambda_{1}+\lambda_{2} \quad \text{and} \quad ad-bc=\lambda_{1}\lambda_{2} $$

then the equation (2.3) can be re-written as,

$$ {}^{J}D^{2\alpha}\left[y\right]-(\lambda_{1}+\lambda_{2})\left({}^{J}D^{\alpha}\left[y\right]\right)+\lambda_{1}\lambda_{2}y=0 \qquad (2.4) $$

For $\lambda_1, \lambda_2$ real and distinct, solution of the equation (2.4) can be written in the form (from Theorem 1)

$$y = A_1 E_\alpha(\lambda_1 t^\alpha) + B_1 E_\alpha(\lambda_2 t^\alpha) \qquad (2.5)$$

Again from second equation of (2.1) we get after putting the value of $y$ the following,

$$
\begin{aligned}
&{}^J D^\alpha [y] = cx + dy \\
&{}^J D^\alpha \left[ A_1 E_\alpha(\lambda_1 t^\alpha) + B_1 E_\alpha(\lambda_2 t^\alpha) \right] = cx + d \left[ A_1 E_\alpha(\lambda_1 t^\alpha) + B_1 E_\alpha(\lambda_2 t^\alpha) \right] \\
&cx = \left[ A_1 \lambda_1 E_\alpha(\lambda_1 t^\alpha) + B_1 \lambda_2 E_\alpha(\lambda_2 t^\alpha \right] - d \left[ A_1 E_\alpha(\lambda_1 t^\alpha) + B_1 E_\alpha(\lambda_2 t^\alpha) \right] \\
&cx = A_1(\lambda_1 - d) E_\alpha(\lambda_1 t^\alpha) + B_1(\lambda_2 - d) E_\alpha(\lambda_2 t^\alpha)
\end{aligned}
$$

From above we obtain the following

$$
\begin{aligned}
x &= \tfrac{1}{c} \left[ A_1(\lambda_1 - d) E_\alpha(\lambda_1 t^\alpha) + B_1(\lambda_2 - d) E_\alpha(\lambda_2 t^\alpha) \right] \\
&= A_2 E_\alpha(\lambda_1 t^\alpha) + B_2 E_\alpha(\lambda_2 t^\alpha) \qquad (2.6)
\end{aligned}
$$

$$A_2 = \frac{A_1(\lambda_1 - d)}{c}; \qquad B_2 = \frac{B_1(\lambda_2 - d)}{c}$$

$A_1, B_1$ are arbitrary constants in above derivation, and $\lambda_1 + \lambda_2 = a + d$, $\lambda_1 \lambda_2 = ad - bc$.

Thus the solution can be written in the form

$$\begin{bmatrix} x \\ y \end{bmatrix} = \begin{bmatrix} A_2 \\ A_1 \end{bmatrix} E_\alpha(\lambda_1 t^\alpha) + \begin{bmatrix} B_2 \\ B_1 \end{bmatrix} E_\alpha(\lambda_2 t^\alpha) \qquad (2.7)$$

Again to solve the system of fractional differential equation (2.1) we use the method similar to as used in classical differential equations.

Since ${}^J D^\alpha [x(t)] = \lambda x(t)$ has solution in the form $x(t) = A E_\alpha(\lambda t^\alpha)$, $A$ is arbitary constant. [15],

putting

$$x = A E_\alpha(\lambda t^\alpha) \quad \text{and} \quad y = B E_\alpha(\lambda t^\alpha)$$

in (2.1) we get the following

$$
\begin{aligned}
&{}^J D^\alpha x - ax - by = 0 \\
&{}^J D^\alpha \left[ A E_\alpha(\lambda t^\alpha) \right] - aA E_\alpha(\lambda t^\alpha) - bB E_\alpha(\lambda t^\alpha) = 0 \\
&A \lambda E_\alpha(\lambda t^\alpha) - A a E_\alpha(\lambda t^\alpha) - bB E_\alpha(\lambda t^\alpha) = 0 \\
&A(\lambda - a) - Bb = 0
\end{aligned}
$$

$$ {}^{J}D^{\alpha} y - cx - dy = 0 $$
$$ {}^{J}D^{\alpha}\left[BE_{\alpha}(\lambda t^{\alpha})\right] - cAE_{\alpha}(\lambda t^{\alpha}) - dBE_{\alpha}(\lambda t^{\alpha}) = 0 $$
$$ B\lambda E_{\alpha}(\lambda t^{\alpha}) - AcE_{\alpha}(\lambda t^{\alpha}) - dBE_{\alpha}(\lambda t^{\alpha}) = 0 $$
$$ -Ac + B(\lambda - d) = 0 $$

$$ \left.\begin{array}{l} A(\lambda - a) - bB = 0 \\ -cA + (\lambda - d)B = 0 \end{array}\right\} \qquad (2.8) $$

Eliminating *A* and *B* from (2.8) we get,

$$ \begin{vmatrix} a-\lambda & b \\ c & d-\lambda \end{vmatrix} = 0 \quad \text{giving} \quad \lambda^2 - (a+d)\lambda + (ad - bc) = 0 \qquad (2.9) $$

which is known as the characteristic equation with roots $\lambda_1$ and $\lambda_2$ , also termed as eigen-values. Three cases may arises

i) The roots $\lambda_1$ and $\lambda_2$ are real and distinct.
ii) The roots are real and equal i.e. $\lambda_1 = \lambda_2 = \lambda$ (say) .
iii) The roots are complex i.e. of the form $\lambda_1, \lambda_2 = p \pm iq$ (say) .

**Case –I**

For real and distinct roots $\lambda_1$ and $\lambda_2$ , we write

$$ x_1 = A_1 E_{\alpha}(\lambda_1 t^{\alpha}),\ y_1 = B_1 E_{\alpha}(\lambda_1 t^{\alpha}) \text{ and } x_2 = A_2 E_{\alpha}(\lambda_2 t^{\alpha}), y_2 = B_2 E_{\alpha}(\lambda_2 t^{\alpha}) $$

The solution of (2.1) is

$$ x = x_1 + x_2 = A_1 E_{\alpha}(\lambda_1 t^{\alpha}) + A_2 E_{\alpha}(\lambda_2 t^{\alpha}) $$
$$ y = y_1 + y_2 = B_1 E_{\alpha}(\lambda_1 t^{\alpha}) + B_2 E_{\alpha}(\lambda_2 t^{\alpha}) $$

Thus the solution of the system of fractional differential equation can be written as,

$$ \begin{bmatrix} x \\ y \end{bmatrix} = c_1 \underset{\sim}{A} E_{\alpha}(\lambda_1 t^{\alpha}) + c_2 \underset{\sim}{B} E_{\alpha}(\lambda_2 t^{\alpha}) $$

Where $\underset{\sim}{A} = \begin{bmatrix} A_1 \\ B_1 \end{bmatrix}$ and $\underset{\sim}{B} = \begin{bmatrix} A_2 \\ B_2 \end{bmatrix}$

are the arbitry constants.

**Example: 1**

$$\left.\begin{array}{l} {}^{J}D^{\alpha}x = 2x + y \\ {}^{J}D^{\alpha}y = x + 2y \end{array}\right\}, \qquad 0 < \alpha \le 1 \qquad \text{with} \qquad x(0) = 2, \qquad y(0) = 1$$

Let

$$x = AE_{\alpha}(\lambda t^{\alpha}) \qquad \text{and} \qquad y = BE_{\alpha}(\lambda t^{\alpha})$$

be the solution of the above differential equation. Substituting this in the above equations we get,

$$\left.\begin{array}{l} A(2-\lambda) + B = 0 \\ A + (2-\lambda)B = 0 \end{array}\right\} \text{......................................................}(2.10)$$

The corresponding characteristic equation is,

$$\begin{vmatrix} 2-\lambda & 1 \\ 1 & 2-\lambda \end{vmatrix} = 0 \qquad \text{giving} \qquad \lambda = 1, 3.$$

Putting $\lambda = 1$ in (2.8) we get $A + B = 0$, taking $A = 1$ we get $B = -1$ and putting $\lambda = 3$ in (2.8) we get $A = B$, taking $A = 1$ we get $B = 1$. Hence the solutions are,

$$x_1 = E_{\alpha}(t^{\alpha}), y_1 = -E_{\alpha}(t^{\alpha}) \ \text{ and } x_2 = E_{\alpha}(3t^{\alpha}), y_2 = E_{\alpha}(3t^{\alpha}).$$

Thus the general solution is

$$\begin{array}{l} x = c_1 E_{\alpha}(t^{\alpha}) + c_2 E_{\alpha}(3t^{\alpha}) \\ y = -c_1 E_{\alpha}(t^{\alpha}) + c_2 E_{\alpha}(3t^{\alpha}) \end{array}.$$

Where $c_1$, $c_2$ are arbitrary constants.

Using the initial condition

$$x(0) = 2, \ y(0) = 0$$

we get $c_1 = c_2 = 1$

Thus the required solution is,

$$\begin{array}{l} x = E_{\alpha}(t^{\alpha}) + E_{\alpha}(3t^{\alpha}) \\ y = -E_{\alpha}(t^{\alpha}) + E_{\alpha}(3t^{\alpha}) \end{array}$$

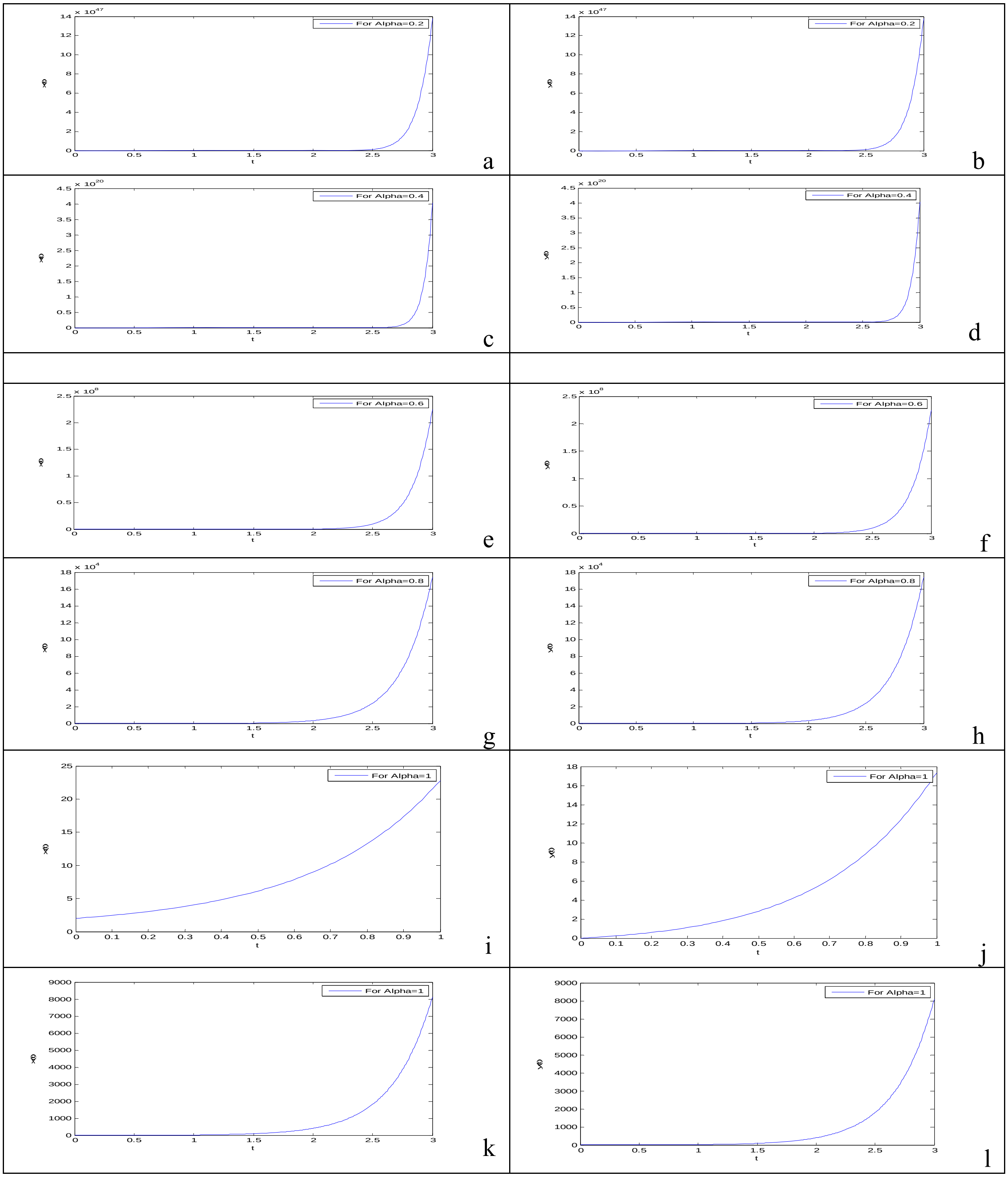

**Fig:1 Numerical simulation of the solutions of fractional differential equation in Example-1 for $x(t)$ and $y(t)$ for different values of $\alpha = 0.2, 0.4, 0.6, 0.8$ and $1.0$. Figure (i) & (k) and figure (j) & (l) represents the same figure only length of x-axis is changed here to represent the prominent initial values.**

Figure-1 represents the graphical presentation of $x(t)$ and $y(t)$ when the eigen-values of the system of differential equations are positive. Numerical simulation shows that $x(t)$ and $y(t)$ both grow rapidly with decrease of order of derivative i.e. as $\alpha$ decreases from 1 towards 0.

**Example: 2**

$$\left.\begin{aligned} {}^{J}D^{\alpha}x &= -2x + y \\ {}^{J}D^{\alpha}y &= x - 2y \end{aligned}\right\}, \qquad 0 < \alpha \le 1 \qquad \text{with} \qquad x(0) = 2, \qquad y(0) = 1$$

Let $\quad x = AE_{\alpha}(\lambda t^{\alpha}) \quad$ and $\quad y = BE_{\alpha}(\lambda t^{\alpha})$

be the solution of the above differential equation. Substituting this in the above equations we get,

$$\left.\begin{aligned} A(-2-\lambda) + B &= 0 \\ A + (-2-\lambda)B &= 0 \end{aligned}\right\} \text{..........................................................(2.10)}$$

The corresponding characteristic equation is,

$$\begin{vmatrix} -2-\lambda & 1 \\ 1 & -2-\lambda \end{vmatrix} = 0 \qquad \text{giving} \qquad \lambda = -1, -3 .$$

Putting $\lambda = -1$ in (2.8) we get $A - B = 0$, taking $A = 1$ we get $B = 1$ and putting $\lambda = -3$ in (2.8) we get $A = -B$ , taking $A = 1$ we get $B = -1$. Hence the solutions are,

$$x_1 = E_{\alpha}(-t^{\alpha}), y_1 = E_{\alpha}(-t^{\alpha}) \text{ and } x_2 = E_{\alpha}(-3t^{\alpha}), y_2 = -E_{\alpha}(-3t^{\alpha})$$

Thus the general solution is

$$\begin{aligned} x &= c_1 E_{\alpha}(-t^{\alpha}) + c_2 E_{\alpha}(-3t^{\alpha}) \\ y &= c_1 E_{\alpha}(-t^{\alpha}) - c_2 E_{\alpha}(-3t^{\alpha}) \cdot \end{aligned}$$

Where $c_1$, $c_2$ are arbitrary constants.

Using the initial condition $\quad x(0) = 2$ , $y(0) = 0$

we get $\quad c_1 = c_2 = 1$ .

Thus the required solution is,

$$\begin{aligned} x &= E_{\alpha}(-t^{\alpha}) + E_{\alpha}(-3t^{\alpha}) \\ y &= E_{\alpha}(-t^{\alpha}) - E_{\alpha}(-3t^{\alpha}) \end{aligned}$$

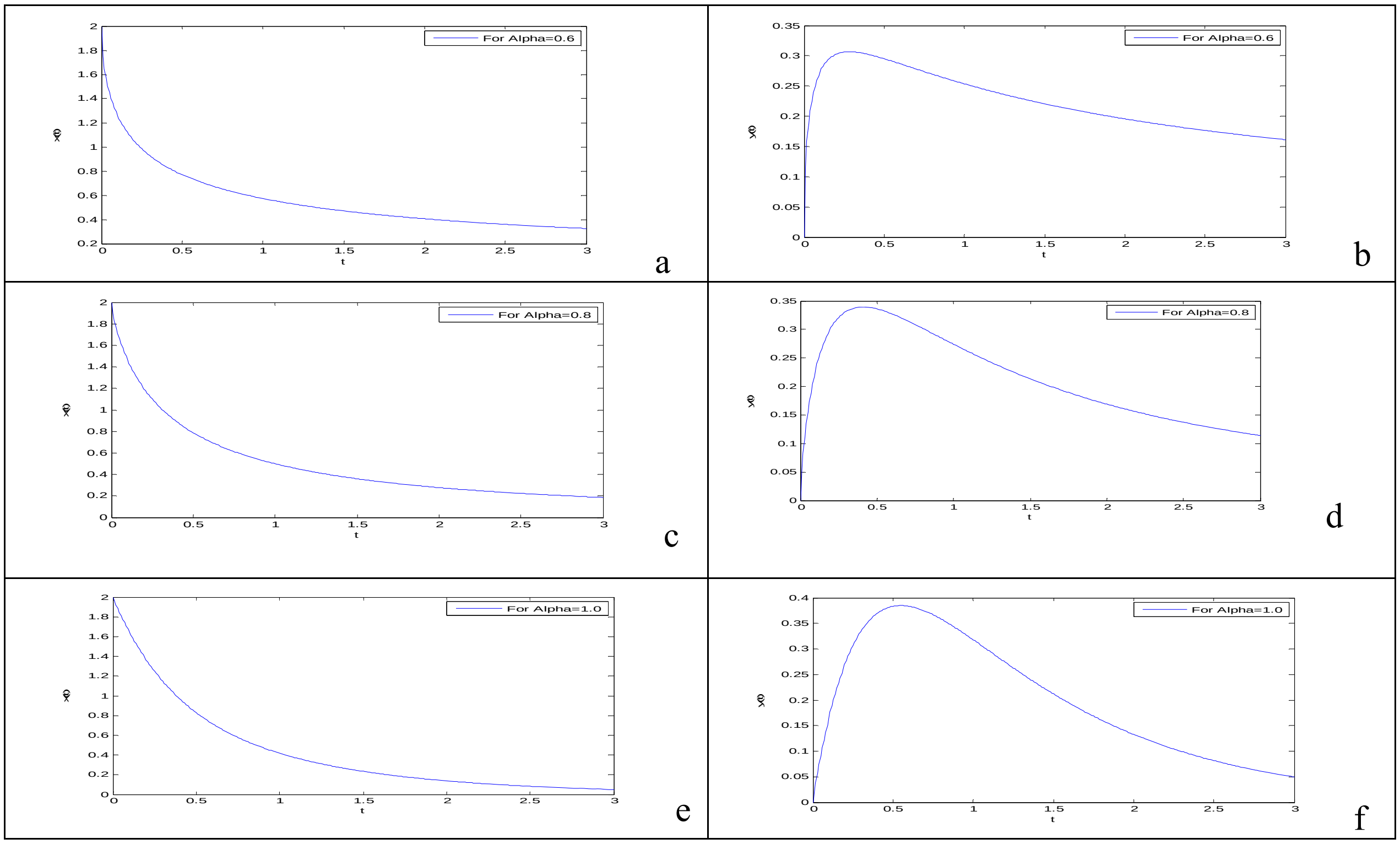

**Fig:2 Numerical simulation of the solutions of fractional differential equation in Example-2 for** $x(t)$ **and** $y(t)$ **for different Values of** $\alpha = 0.6, 0.8$ and $1.0$.

Figure-2 represents the graphical presentation of $x(t)$ and $y(t)$ when the eigen-values of the system of differential equations are negative. Numerical simulation shows that $x(t)$ and $y(t)$ are both decaying rapidly with decrease of order of derivative i.e. as $\alpha$ decreases from 1 to 0.6. For negative eigen-values the solutions are decaying asymptotically to zero.

## Case –II

The roots of the equation (2.3) are complex and are of the form $\lambda_1, \lambda_2 = p \pm iq$ then the solution $x(t)$ and $y(t)$ can be written in the form (from Theorem 3),

$$x = (A_1 + iA_2)\left[ E_\alpha (p+iq)t^\alpha \right] \quad \text{and} \quad y = (B_1 + iB_2)\left[ E_\alpha (p+iq)t^\alpha \right]$$

Using the definition of Mittag-Leffler function and fractional cosine and sine functions, that is

$$E_\alpha(ix^\alpha) \overset{\text{def}}{=} \cos_\alpha(x^\alpha) + i\sin(x^\alpha)$$

we get,

$$x=(A_1+iA_2)\left[E_\alpha(p+iq)t^\alpha\right]=(A_1+iA_2)\left[E_\alpha(pt^\alpha)\right]\left[E_\alpha(iq)t^\alpha\right]$$
$$=(A_1+iA_2)\left[E_\alpha(pt^\alpha)\right]\left[\cos_\alpha(qt^\alpha)+i\sin_\alpha(qt^\alpha)\right]$$
$$=\left[E_\alpha(pt^\alpha)\right]\left[\left(A_1\cos_\alpha(qt^\alpha)-A_2\sin_\alpha(qt^\alpha)\right)+i\left(A_2\cos_\alpha(qt^\alpha)+A_1\sin_\alpha(qt^\alpha)\right)\right]$$

Similarly we get by repeating the above steps for $y$ as follows

$$y=\left[E_\alpha(pt^\alpha)\right]\left[\left(B_1\cos_\alpha(qt^\alpha)-B_2\sin_\alpha(qt^\alpha)\right)+i\left(B_2\cos_\alpha(qt^\alpha)+B_1\sin_\alpha(qt^\alpha)\right)\right]$$

In above obtained expressions for $x$ and $y$, we have complex quantity as $u+iv$. This may be also considered as linear combination of $u$ and $v$ considering $i$ a constant. Therefore, we can say that $x$ is linear combination of $x_1$ (the real part of obtained complex $x$), and $x_2$ (the imaginary part of obtained complex $x$). Similarly we have $y$ as linear combination of $y_1$ and $y_2$ [21]. With this argument we write the following

$$x_1=\left[E_\alpha(pt^\alpha)\right]\left[A_1\cos_\alpha(qt^\alpha)-A_2\sin_\alpha(qt^\alpha)\right]$$
$$x_2=\left[E_\alpha(pt^\alpha)\right]\left[A_2\cos_\alpha(qt^\alpha)+A_1\sin_\alpha(qt^\alpha)\right]$$

$$y_1=\left[E_\alpha(pt^\alpha)\right]\left[B_1\cos_\alpha(qt^\alpha)-B_2\sin_\alpha(qt^\alpha)\right]$$
$$y_2=\left[E_\alpha(pt^\alpha)\right]\left[B_2\cos_\alpha(qt^\alpha)+B_1\sin_\alpha(qt^\alpha)\right]$$

It can be shown that $(x_1,y_1);(x_2,y_2)$ are solutions of the given equations (2.1).
Thus the general solution in this case can be written in the form as in classical integer order differential equation [21, pp.305].

The linear combination of $x_1$ and $x_2$, gives $x$ and linear combination of $y_1$ and $y_2$ gives $y$, which is represented as following

$$x=\left[E_\alpha(pt^\alpha)\right]\left[M\left(A_1\cos_\alpha(qt^\alpha)-A_2\sin_\alpha(qt^\alpha)\right)+N\left(A_2\cos_\alpha(qt^\alpha)+A_1\sin_\alpha(qt^\alpha)\right)\right]$$

$$y=\left[E_\alpha(pt^\alpha)\right]\left[M\left(B_1\cos_\alpha(qt^\alpha)-B_2\sin_\alpha(qt^\alpha)\right)+N\left(B_2\cos_\alpha(qt^\alpha)+B_1\sin_\alpha(qt^\alpha)\right)\right]$$

With $M,N$ as arbitrary constants, determined from initial states. We demonstrate by following examples.

**Example: 3**

$$\left.\begin{array}{l}{}^JD^\alpha x=3x+2y\\ {}^JD^\alpha y=-5x+y\end{array}\right\}\qquad 0<\alpha\le 1\qquad \text{with}\qquad x(0)=2,\qquad y(0)=1$$

Let

$$x=AE_\alpha(\lambda t^\alpha)\qquad\text{and}\qquad y=BE_\alpha(\lambda t^\alpha)$$

be the solution of the above differential equation, then putting in the above equation we get,

$$\left.\begin{aligned} A(3-\lambda)+2B=0 \\ -5A+(1-\lambda)B=0 \end{aligned}\right\} \qquad \text{..........}(2.11)$$

The corresponding characteristic equation is,

$$\begin{vmatrix} 3-\lambda & 2 \\ -5 & 1-\lambda \end{vmatrix}=0 \quad \text{giving} \quad \lambda=2\pm 3i\,.$$

Putting $\lambda=2+3i$ in (2.8) by putting $a=3, b=2, c=-5, d=1$ we get the following

$$\begin{aligned}
&A(\lambda-a)-2B=0\\
&A(2+3i-3)+2B=0\\
&2B=(-1+3i)A\\
&-cA+(\lambda-d)B=0\\
&5A+(2+3i-1)B=0\\
&5A=-(1+3i)B\\
&5A\times(1-3i)=-(1+3i)\times(1-3i)B\\
&5(1-3i)A=-10B\\
&2B=(-1+3i)A
\end{aligned}$$

The (2.8) returns the same answer that is $2B=(3i-1)A$.We choose here $A=2$, so $B=3i-1$, as we obtained one equation with two unknowns. Thus $B=3i-1$ and $A=2$ can be taken as one of the trial solution of the above. Hence the solution is,

$$\begin{aligned}
x&=2E_{\alpha}\left((2+3i)t^{\alpha}\right)\\
&=2\left[E_{\alpha}(2t^{\alpha})\right]\left[\cos_{\alpha}(3t^{\alpha})+i\sin_{\alpha}(3t^{\alpha})\right]
\end{aligned}$$

Thus

$$x_1=2\left[E_{\alpha}(2t^{\alpha})\right]\left[\cos_{\alpha}(3t^{\alpha})\right] \quad \text{and} \quad x_2=2\left[E_{\alpha}(2t^{\alpha})\right]\left[\sin_{\alpha}(3t^{\alpha})\right]$$

Similarly the solution for $y$ can be written in the form

$$\begin{aligned}
y&=(-1+3i)\left[E_{\alpha}(2+3i)t^{\alpha}\right]\\
&=(-1+3i)\left[E_{\alpha}(2t^{\alpha})\right]\left[\cos_{\alpha}(3t^{\alpha})+i\sin_{\alpha}(3t^{\alpha})\right]\\
&=\left[E_{\alpha}(2t^{\alpha})\right]\left[\left(-\cos_{\alpha}(3t^{\alpha})-3\sin_{\alpha}(3t^{\alpha})\right)+i\left(3\cos_{\alpha}(3t^{\alpha})-\sin_{\alpha}(3t^{\alpha})\right)\right]
\end{aligned}$$

Hence $y_1=\left[E_{\alpha}(2t^{\alpha})\right]\left[-\cos_{\alpha}(3t^{\alpha})-3\sin_{\alpha}(3t^{\alpha})\right]$ and $y_2=\left[E_{\alpha}(2t^{\alpha})\right]\left[3\cos_{\alpha}(3t^{\alpha})-\sin_{\alpha}(3t^{\alpha})\right]$

Therefore the general solution is linear combination of $x_1, x_2$ for $x$ and linear combination of $y_1$, $y_2$ for $y$, and we write the following

$$x = \left[E_\alpha(2t^\alpha)\right]\left[2M\cos_\alpha(3t^\alpha) + 2N\sin_\alpha(3t^\alpha)\right]$$
$$y = \left[E_\alpha(2t^\alpha)\right]\left[M\left(-\cos_\alpha(3t^\alpha) - 3\sin_\alpha(3t^\alpha)\right) + N\left(3\cos_\alpha(3t^\alpha) - \sin_\alpha(3t^\alpha)\right)\right]$$

where $M$ and $N$ are arbitrary constants. Using initial conditions $x(0) = 2, y(0) = 1$ we get $2M = 2,$ and $3N - M = 1$ giving $M = 1, N = 2$.

Hence the required solution is,

$$x = \left[E_\alpha(2t^\alpha)\right]\left[2\cos_\alpha(3t^\alpha) + 4\sin_\alpha(3t^\alpha)\right]$$
$$y = \left[E_\alpha(2t^\alpha)\right]\left[\cos_\alpha(3t^\alpha) - 5\sin_\alpha(3t^\alpha)\right]$$

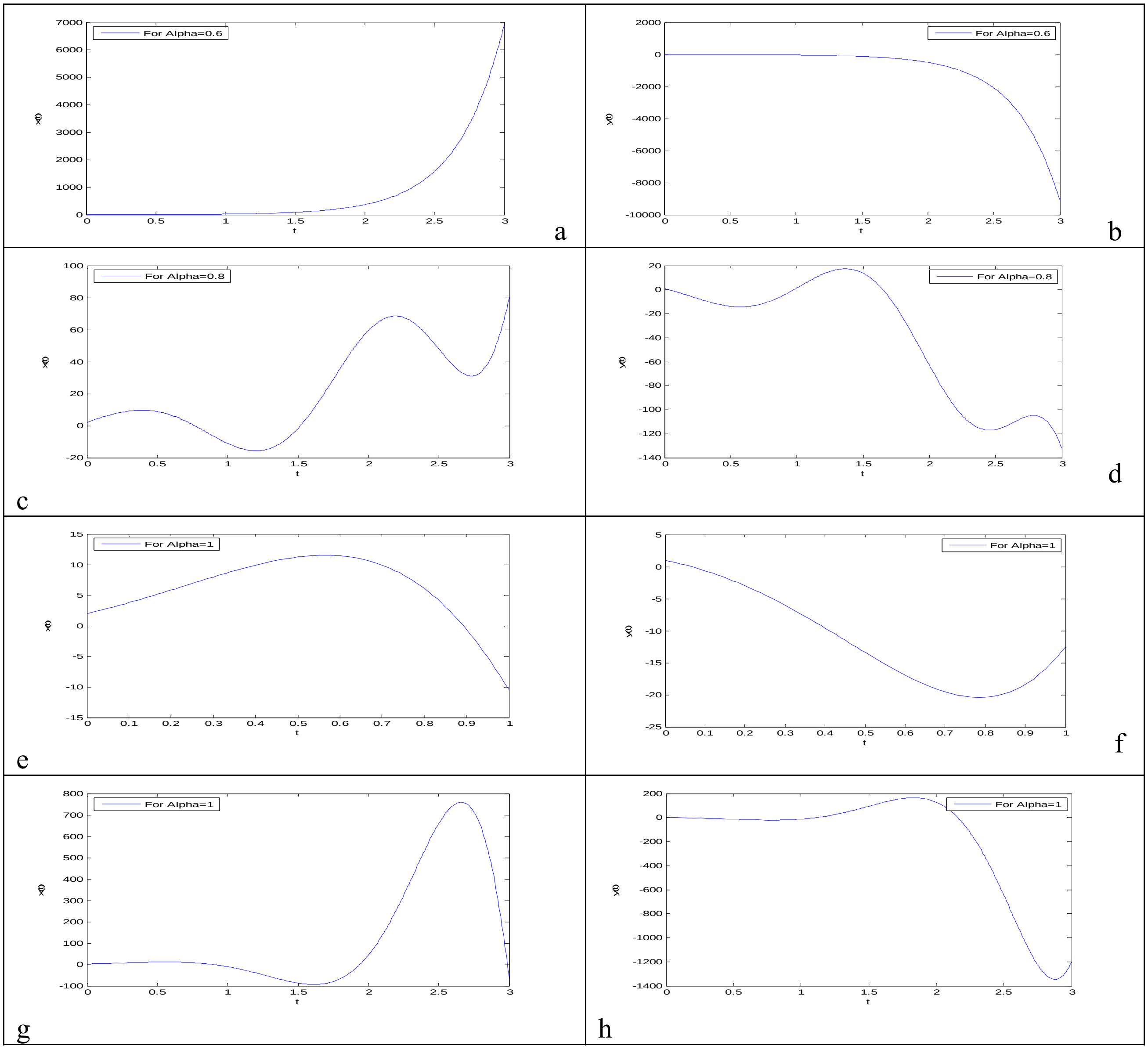

**Fig:3 Numerical simulation of the solutions of fractional differential equation in Example-3 for $x(t)$ and $y(t)$ for different Values of $\alpha = 0.6, 0.8$ and $1.0$. Figure (e) & (f) and Figure (g) & (h) represents the same figure only length of x-axis is changed here to represent the prominent initial values.**

Numerical simulation in figure-3 shows that for $\alpha = 0.6$, $0.8$ and $1.0$ after the initiation $(t = 0)$ of the system $x(t)$ and $y(t)$ both oscillate, Period of oscillation changes with decrease of $\alpha$.

## Case-III

In this case roots of the equations (2.9) being equal, that is $\lambda_1 = \lambda_2 = \lambda$ . Then one solution will be of the form

$$x_1 = AE_\alpha(\lambda t^\alpha) \quad \text{and} \quad y_1 = BE_\alpha(\lambda t^\alpha)$$

and the other solution will be

$$x_2 = (A_1 t^{\alpha} + A_2) E_{\alpha}(\lambda t^{\alpha}) \quad \text{and} \quad y_2 = (B_1 t^{\alpha} + B_2) E_{\alpha}(\lambda t^{\alpha})$$

Hence the general solution is,

$$x = c_1 A E_{\alpha}(\lambda t^{\alpha}) + c_2 (A_1 t^{\alpha} + A_2) E_{\alpha}(\lambda t^{\alpha})$$
$$y = c_1 B E_{\alpha}(\lambda t^{\alpha}) + c_2 (B_1 t^{\alpha} + B_2) E_{\alpha}(\lambda t^{\alpha})$$

where $A, A_1, A_2, B, B_1, B_2, c_1, c_2$ are arbitry constants

**Example4:**

$$\left.\begin{aligned} {}^{J}D^{\alpha} x &= 4x - y \\ {}^{J}D^{\alpha} y &= x + 2y \end{aligned}\right\} \quad 0 < \alpha \le 1 \quad \text{with} \quad x(0) = 2, \quad y(0) = 1 \text{...........(2.12)}$$

Let

$$x = A E_{\alpha}(\lambda t^{\alpha}) \quad \text{and} \quad y = B E_{\alpha}(\lambda t^{\alpha})$$

be the solutions of the above differential equation, then putting in the above equation we get

$$\left.\begin{aligned} A(4-\lambda) - B &= 0 \\ A + (2-\lambda) B &= 0 \end{aligned}\right\} \text{..........................................................(2.13)}$$

Eliminating $A$ and $B$ as in previous examples, we get,

$$\begin{vmatrix} 4-\lambda & 1 \\ 1 & 2-\lambda \end{vmatrix} = 0 \quad \text{giving} \quad \lambda = 3,3$$

For $\lambda = 3$ from equation (2.4) we get $A = B = 1$.

Thus

$$x_1 = E_{\alpha}(3t^{\alpha}) \quad \text{and} \quad y_1 = E_{\alpha}(3t^{\alpha})$$

is one solution of the equation. The second solution is as in classical integer order differential equation [21, pp.307]

$$x_2 = (A_1 t^{\alpha} + A_2)\left[E_{\alpha}(3t^{\alpha})\right] \quad \text{and} \quad y_2 = (B_1 t^{\alpha} + B_2)\left[E_{\alpha}(3t^{\alpha})\right]$$

$\alpha$ -th order differentiating for above $x$ and $y$ we obtain the following

$${}^{J}D^{\alpha} x = 3(A_1 t^{\alpha} + A_2) E_{\alpha}(3t^{\alpha}) + \Gamma(1+\alpha) A_1 E_{\alpha}(3t^{\alpha})$$

and

$$^{J}D^{\alpha} y = 3(B_1 t^{\alpha} + B_2) E_{\alpha}(3t^{\alpha}) + \Gamma(1+\alpha) B_1 E_{\alpha}(3t^{\alpha})$$

Putting the above obtained result in the given equation (2.12) we get,

$$3(A_1 t^{\alpha} + A_2) E_{\alpha}(3t^{\alpha}) + \Gamma(1+\alpha) A_1 E_{\alpha}(3t^{\alpha}) = E_{\alpha}(3t^{\alpha})(4A_1 - B_1) t^{\alpha} + E_{\alpha}(3t^{\alpha})(4A_2 - B_2)$$
$$3(B_1 t^{\alpha} + B_2) E_{\alpha}(3t^{\alpha}) + \Gamma(1+\alpha) B_1 E_{\alpha}(3t^{\alpha}) = E_{\alpha}(3t^{\alpha})(A_1 + 2B_1) t^{\alpha} + E_{\alpha}(3t^{\alpha})(A_2 + 2B_2)$$

$$3(A_1 t^{\alpha} + A_2) + \Gamma(1+\alpha) A_1 = (4A_1 - B_1) t^{\alpha} + (4A_2 - B_2)$$
$$3(B_1 t^{\alpha} + B_2) + \Gamma(1+\alpha) B_1 = (A_1 + 2B_1) t^{\alpha} + (A_2 + 2B_2)$$

Comparing the coefficients and simplifying we get

$$A_1 = B_1 \qquad A_2 - B_2 = A_1 \Gamma(1+\alpha) = B_1 \Gamma(1+\alpha)$$

for simple non-zero values we take

$A_1 = B_1 = 1, \qquad A_2 - B_2 = \Gamma(1+\alpha),$ For simplicity take $B_2 = 0$ then $A_2 = \Gamma(1+\alpha)$.

Thus the other solution is

$$x = \left[t^{\alpha} + \Gamma(1+\alpha)\right]\left[E_{\alpha}(3t^{\alpha})\right] \quad \text{and} \quad y = t^{\alpha}\left[E_{\alpha}(3t^{\alpha})\right]$$

Hence the general solution can be written in the form as in classical integer order differential equation [21, pp.307]

$$x = c_1 E_{\alpha}(3t^{\alpha}) + c_2\left[t^{\alpha} + \Gamma(1+\alpha)\right] E_{\alpha}(3t^{\alpha}) \quad \text{and} \quad y = c_1 E_{\alpha}(3t^{\alpha}) + c_2 t^{\alpha} E_{\alpha}(3t^{\alpha})$$

Putting the initial condition.

$$x(0) = 2 \quad \text{and} \quad y(0) = 1$$

and solving we get

$$c_1 = 1 \quad \text{and} \quad c_2 = \frac{1}{\Gamma(1+\alpha)}$$

Hence the solution is,

$$x = E_{\alpha}(3t^{\alpha}) + \tfrac{1}{\Gamma(1+\alpha)}\left[t^{\alpha} + \Gamma(1+\alpha)\right]\left[E_{\alpha}(3t^{\alpha})\right] \quad \text{and} \quad y = E_{\alpha}(3t^{\alpha}) + \tfrac{1}{\Gamma(1+\alpha)} t^{\alpha} E_{\alpha}(3t^{\alpha})$$

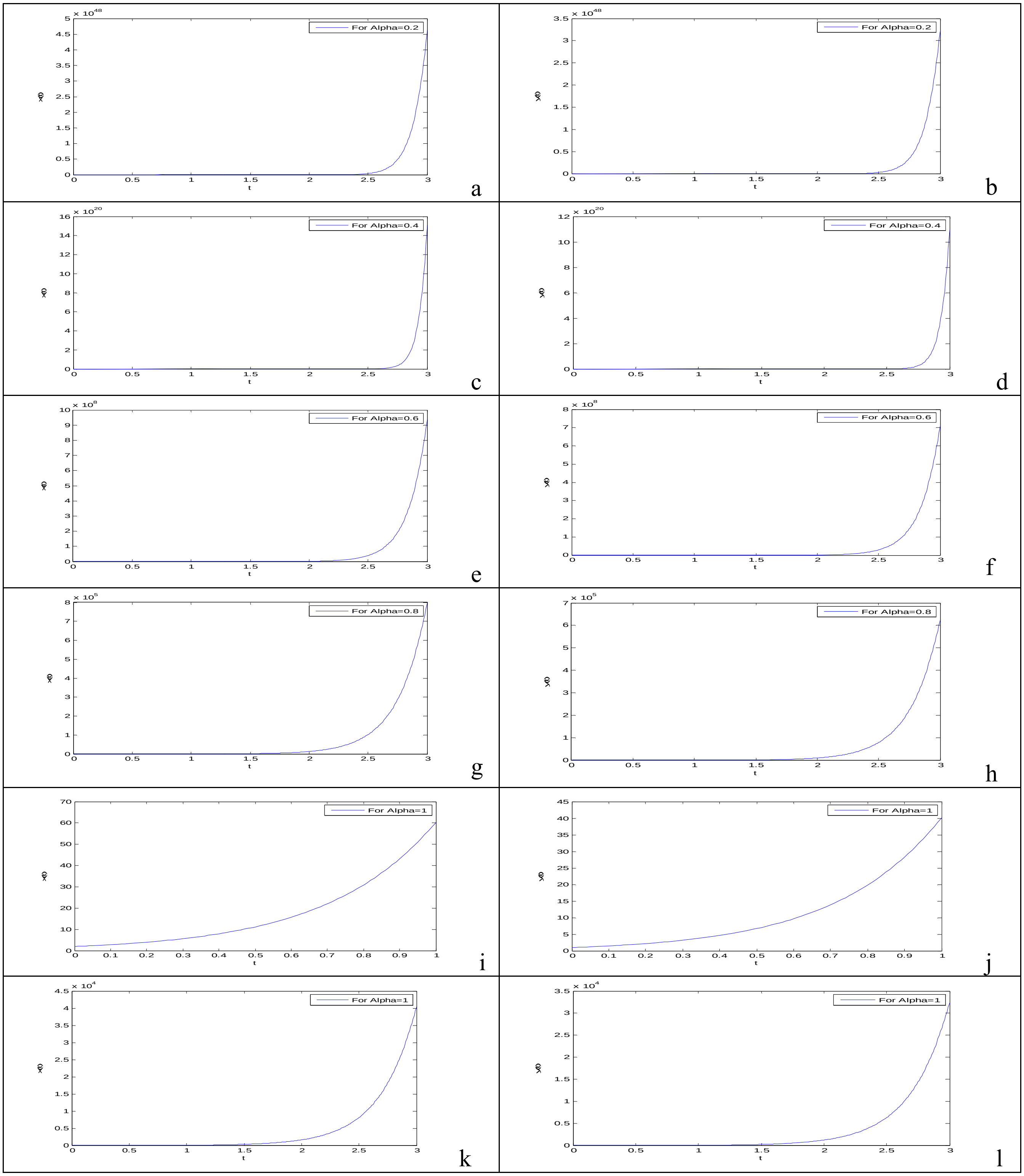

**Fig:4 Numerical simulation of the solutions of fractional differential equation in Example-3 for $x(t)$ and $y(t)$ for different Values of $\alpha = 0.2, 0.4, 0.6, 0.8 \text{ and } 1.0.$ Figure (i) & (k) and Figure (j) & (l) represents the same figure only length of x-axis is changed here to represent the prominent initial values.**

Numerical simulation in figure-4 shows that $x(t)$ and $y(t)$ start from 2 and 1 respectively; and both of them start to grow. This time interval to grow decreases as $\alpha$ increases. Moreover growth of $x(t)$ and $y(t)$ is higher for is for lower $\alpha$ values. As $\alpha$ increases growth rate of the solution decreases. This implies for lower values of $\alpha$, $x(t)$ and $y(t)$ grow initially slowly. Once they start to grow, their growth rate is very high whereas for higher $\alpha$ values $x(t)$ and $y(t)$ start to grow sooner but their growth rate is low.

From the above discussion of the three cases one can state a theorem in the following form

**Theorem 4**:

The solutions of the system of differential equations

$$^{J}D^{\alpha}X = AX, \qquad \text{Where} \qquad A = \begin{bmatrix} a & b \\ c & d \end{bmatrix}, \qquad X = \begin{bmatrix} x \\ y \end{bmatrix}$$

are

Case (i)

$$X = c_1 \underset{\sim}{A} E_{\alpha}(\lambda_1 t^{\alpha}) + c_2 \underset{\sim}{B} E_{\alpha}(\lambda_2 t^{\alpha}),$$

Where $\underset{\sim}{A}, \underset{\sim}{B}, c_1$ and $c_2$ are the arbitry constants,

$$\lambda_1 + \lambda_2 = a + d,$$
$$\lambda_1 \lambda_2 = ad - bc,$$
$$\lambda_1 \neq \lambda_2.$$

Case (ii)

$$X = c_1 \underset{\sim}{A} E_{\alpha}(\lambda t^{\alpha}) + c_2 (\underset{\sim}{B}_1 t^{\alpha} + \underset{\sim}{B}_2) E_{\alpha}(\lambda t^{\alpha}),$$

Where $\underset{\sim}{A}, \underset{\sim}{B}_1, \underset{\sim}{B}_2, c_1$ and $c_2$ are the arbitry constants,

$$2\lambda = a + d,$$
$$\lambda^2 = ad - bc.$$

Case (iii)

$$X = \left[ E_{\alpha}(pt^{\alpha}) \right] \begin{bmatrix} A_1 \cos_{\alpha}(qt^{\alpha}) - A_2 \sin_{\alpha}(qt^{\alpha}) + (A_2 \cos_{\alpha}(qt^{\alpha}) + A_1 \sin_{\alpha}(qt^{\alpha})) \\ B_1 \cos_{\alpha}(qt^{\alpha}) - B_2 \sin_{\alpha}(qt^{\alpha}) + (B_2 \cos_{\alpha}(qt^{\alpha}) + B_1 \sin_{\alpha}(qt^{\alpha})) \end{bmatrix}$$

Where $A_1, A_2, B_1, B_2$ are the arbitry constants,

$$2p = a + d,$$
$$p^2 + q^2 = ad - bc.$$

## 3.0 Application of the above formulation in real life problem:

Consider the following fractional damped oscillator, formulated by Jumarie fractional derivative

$$ {}^{J}D^{2\alpha}\left[x\right]+2a\left({}^{J}D^{\alpha}\left[x\right]\right)+bx=0 \qquad (3.1)$$

Let ${}^{J}D^{\alpha}\left[x\right]\equiv \frac{d^{\alpha}x}{dt^{\alpha}}=y$ then the given equation reduce to the following system of equation

$$\left.\begin{aligned} {}^{J}D^{\alpha}\left[y\right]&=\frac{d^{\alpha}y}{dt^{\alpha}}=-2ay-bx \\ {}^{J}D^{\alpha}\left[x\right]&=\frac{d^{\alpha}x}{dt^{\alpha}}=y \end{aligned}\right\} \qquad (3.2)$$

The above system of equation can be written in the form

$$\begin{bmatrix} {}^{J}D^{\alpha}x \\ {}^{J}D^{\alpha}y \end{bmatrix}=\begin{bmatrix} \frac{d^{\alpha}x}{dt^{\alpha}} \\ \frac{d^{\alpha}y}{dt^{\alpha}} \end{bmatrix}=\begin{bmatrix} 0 & 1 \\ -b & -2a \end{bmatrix}\begin{bmatrix} x \\ y \end{bmatrix}$$

Let

$$x=AE_{\alpha}(\lambda t^{\alpha}) \quad \text{and} \quad y=BE_{\alpha}(\lambda t^{\alpha})$$

be solutions of the differential equations.

Then

$$\left.\begin{aligned} A(0-\lambda)+B&=0 \\ -bA+(-2a-\lambda)B&=0 \end{aligned}\right\} \qquad (3.3)$$

For the above system of equation the auxiliary equation is,

$$\begin{vmatrix} 0-\lambda & 1 \\ -b & -2a-\lambda \end{vmatrix}=0 \quad \text{giving} \quad \lambda^{2}+2a\lambda+b=0.$$

Here the discriminant is $4(a^{2}-b)$. We consider the case when $a^{2}-b<0$. then the eigen-values are

$$\lambda_{1},\lambda_{2}=p\pm iq \qquad \text{where} \qquad p=-a, \qquad q=\sqrt{b-a^{2}}\,.$$

Then from (3.3) putting $\lambda=p+iq$ we get $A(p+iq)=B$, we can take the solution in the form $A=1, B=(p+iq)$. The general solution will be of the form,

$$x=E_{\alpha}\left[(p+iq)t^{\alpha}\right]=\left[E_{\alpha}(pt^{\alpha})\right]\left[\cos_{\alpha}(qt^{\alpha})+\sin_{\alpha}(qt^{\alpha})\right]$$

$$x_{1}=\left[E_{\alpha}(pt^{\alpha})\right]\left[\cos_{\alpha}(qt^{\alpha})\right] \quad \text{and} \quad x_{2}=\left[E_{\alpha}(pt^{\alpha})\right]\left[\sin_{\alpha}(qt^{\alpha})\right]$$

Thus the general solution can be written in the form as in classical integer order differential equation [21, pp.305]

$$x=\left[E_{\alpha}(pt^{\alpha})\right]\left[C_{1}\cos_{\alpha}(qt^{\alpha})+C_{2}\sin_{\alpha}(qt^{\alpha})\right]$$

Where $C_1$ and $C_2$ are arbitrary constants.

Here $x(0)=2$ and ${}^{J}D^{\alpha}x(0)=1$ and solving we get $C_{1}=1$ and $C_{2}=(1-2p)/q$. Hence the solution is

$$x=\left[E_{\alpha}(pt^{\alpha})\right]\left[\cos_{\alpha}(qt^{\alpha})+\tfrac{1-2p}{q}\sin_{\alpha}(qt^{\alpha})\right]$$

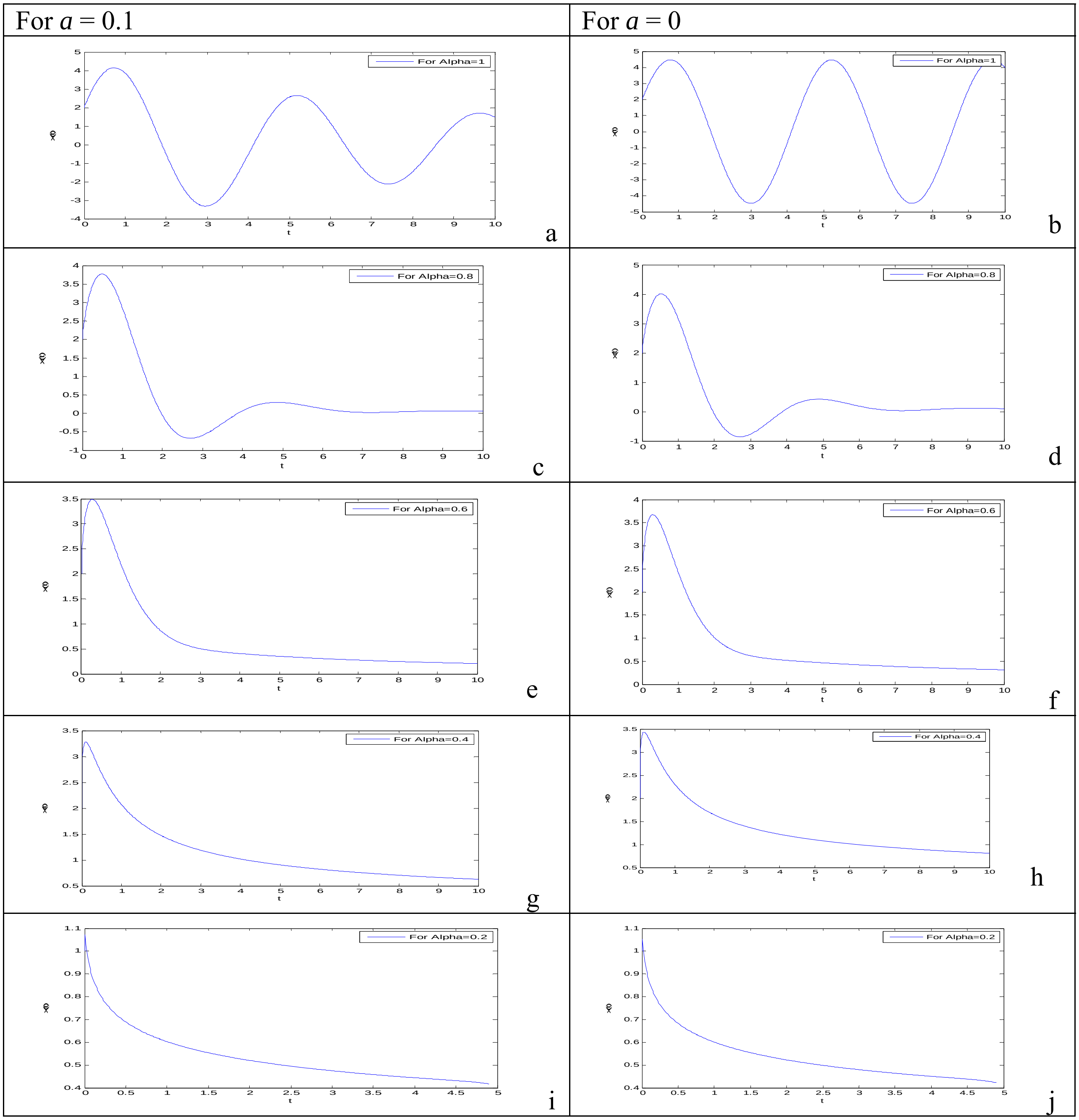

**Fig: 5: Numerical simulation of the solutions of fractional differential equation in equation (3.1) for different Values** $\alpha = 0.2, 0.4, 0.6, 0.8$ and $1.0$. **for *b*=2 and for *a*=0.1 & *a*=0.**

The numerical simulation of the solutions of the differential equation (3.1) has shown in figure 5, for $b = 2$, left hand figures for $a = 0.1$ and right hand figures for $a = 0$. Figures 4 (a) and (b) are drawn for $\alpha = 1$, it is clear from the figure in presence of damping the amplitude of the oscillation decreases with time. Figures 4 (c) and (d) are drawn for $\alpha = 0.8$, it is clear from the figure in both the cases the amplitude of the oscillation decreases with time and ultimately amplitude tends to zero. Figures 4 (e) & (f) and (g) & (h) and (i) & (j) are drawn for $\alpha = 0.6, 0.4$ and $0.2$ respectively, it is clear from the figures with decrease of order of derivative the oscillator losses the oscillating behavior.

## 4.0 Conclusions

The system of fractional differential equation arises in different applications. Here we develop an algorithm to solve the system of fractional differential equations with modified fractional derivative (with Jumarie's fractional derivative formulation) in terms of Mittag-Leffler function and the generalized Sine and Cosine functions. From the numerical simulations it is observed that the growing or decaying of the solutions is fast in fractional order derivative case compare to the integer order derivative. The use of this type of Jumarie fractional derivative gives a conjugation with classical methods of solution of system of linear integer order differential equations, by usage of Mittag-Leffler and generalized trigonometric functions that we have demonstrated here in this paper. The ease of this method and its conjugation to classical method to solve system of linear fractional differential equation is appealing to researchers in fractional dynamic systems.

## Acknowledgement

Acknowledgments are to **Board of Research in Nuclear Science** (BRNS), Department of Atomic Energy Government of India for financial assistance received through BRNS research project no. 37(3)/14/46/2014-BRNS with BSC BRNS, title "Characterization of unreachable (Holderian) functions via Local Fractional Derivative and Deviation Function".